\documentclass[12pt, twoside]{amsart}
\usepackage{parskip}
\usepackage{layout}
\usepackage{amsthm,amsmath,amssymb,oldgerm}
\usepackage[a4paper, margin=2.40cm]{geometry}
\usepackage{mathrsfs}
\usepackage{mathtools}
\usepackage{verbatim} 
\usepackage{todonotes}
\usepackage{hyperref}
\usepackage{bm}

\theoremstyle{plain}
\newtheorem{thm}{Theorem}[section]
\newtheorem*{theorem*}{Main theorem}

\newtheorem{lem}[thm]{Lemma}
\newtheorem{prop}[thm]{Proposition}
\newtheorem{cor}[thm]{Corollary}
\theoremstyle{definition}
\theoremstyle{definition}
\newtheorem{rem}[thm]{Remark}

\title[$L^{\infty}$-norm bounds for Jacobi cusp forms]{$L^{\infty}$-norm bounds for Jacobi cusp forms}
\date{\today}

\begin{document}

\author[A.~Aryasomayajula]{Anilatmaja Aryasomayajula}
\address{Indian Institute of Science Education and Research, Tirupati, Karakambadi Road, Mangalam (P.O.) Tirupati-517507, Andhra Pradesh, India}
\email{anilatmaja@gmail.com}

\author[J.~Kramer]{J\"urg Kramer}
\address{Institut f\"ur Mathematik, Humboldt-Universit\"at zu Berlin, Unter den Linden 6, D-10099 Berlin, Germany}
\email{kramer@math.hu-berlin.de}

\author[A-M.~von Pippich]{Anna-Maria von Pippich}
\address{Fachbereich Mathematik und Statistik, Universit\"at Konstanz, Universit\"atsstra{\ss}e 10, D-78464 Konstanz, Germany}
\email{anna.pippich@uni-konstanz.de}

\begin{abstract}
In this article, we give $L^{\infty}$-norm bounds for the natural invariant norm of cusp forms of real weight $k$ and character $\chi$ for any cofinite 
Fuchsian subgroup $\Gamma\subset\mathrm{SL}_{2}(\mathbb{R})$. Using the representation of Jacobi cusp forms of integral weight $k$ and index 
$m$ for the modular group $\Gamma_{0}=\mathrm{SL}_{2}(\mathbb{Z})$ as linear combinations of modular forms of weight $k-\frac{1}{2}$ for some 
congruence subgroup of $\Gamma_{0}$ (depending on $m$) and suitable Jacobi theta functions, we derive $L^{\infty}$-norm bounds for the natural 
invariant norm of these Jacobi cusp forms. More specifically, letting $J_{k,m}^{\mathrm{cusp}}(\Gamma_{0})$ denote the complex vector space of 
Jacobi cusp forms under consideration and $\Vert\cdot\Vert_{\mathrm{Pet}}$ the pointwise Petersson norm on $J_{k,m}^{\mathrm{cusp}}(\Gamma_
{0})$, we prove that for $k\in\mathbb{Z}_{\ge 5}$ and $m\in\mathbb{Z}_{\ge 1}$, and a given $\epsilon>0$, the $L^{\infty}$-norm bound 
\begin{align*}
\Vert\phi\Vert_{L^{\infty}}=\sup_{(\tau,z)\in\mathbb{H}\times\mathbb{C}}\Vert\phi(\tau,z)\Vert_{\mathrm{Pet}}=O_{\Gamma_{0},\epsilon}\big(k\,m^{\frac
{7}{4}+\epsilon}\big)
\end{align*}
holds for any $\phi\in J_{k,m}^{\mathrm{cusp}}(\Gamma_{0})$, which is $L^{2}$-normalized with respect to the Petersson inner product, where the 
implied constant depends on $\Gamma_{0}$ and the choice of $\epsilon>0$.
\end{abstract}

\maketitle

\section{Introduction}

\label{sec-1}

\subsection{Background} 
\label{subsec1.1} 
In general, bounds for automorphic forms and for their Fourier coefficients represent an area of great interest in number theory. More specifically, we 
mention in this respect the results of~\cite{fjk2}, where J.~Friedman, J.~Jorgenson, and J.~Kramer obtained optimal $L^{\infty}$-norm bounds on average 
for cusp forms of even weight $k$ for any cofinite Fuchsian subgroup $\Gamma\subset\mathrm{SL}_{2}(\mathbb{R})$. These $L^{\infty}$-norm bounds 
turn out to be uniform with respect to the subgroup $\Gamma$. Moreover, in~\cite{fjk1}, effective versions for these $L^{\infty}$-norm bounds are given. 
With regard to $L^{\infty}$-norm bounds for individual Hecke eigenforms of large level, we mention, for example, the results by V.~Blomer and R.~Holowinsky 
in~\cite{blomer-holowinsky}.

So far, less attention has been devoted to the study of $L^{\infty}$-norm bounds for Jacobi cusp forms. The first comprehensive study of Jacobi forms 
was undertaken by M.~Eichler and D.~Zagier in~\cite{eich-zag}. Subsequently, various authors have built on their work. In contrast to their analytical
approach, a geometrical approach to the theory of Jacobi forms was given by J.~Kramer in~\cite{k1}.

Let $k,m$ be positive integers. A Jacobi form of weight $k$ and index $m$ for the modular group $\Gamma_{0}\coloneqq\mathrm{SL}_{2}(\mathbb{Z})$ 
is a holomorphic function on the product $\mathbb{H}\times\mathbb{C}$ of the upper half-plane $\mathbb{H}$ with the complex plane $\mathbb{C}$ 
having a suitable transformation behaviour with respect to $\Gamma_{0}$ and vanishing ``at infinity''. We denote the complex vector space of Jacobi 
cusp forms of weight $k$ and index $m$ for $\Gamma_{0}$ by $J_{k,m}^{\mathrm{cusp}}(\Gamma_{0})$. The pointwise Petersson norm of a Jacobi 
form $\phi\in J_{k,m}^{\mathrm{cusp}}(\Gamma_{0})$ is then defined by
\begin{align*}
\Vert\phi(\tau,z)\Vert^{2}_{\mathrm{Pet}}\coloneqq\vert\phi(\tau,z)\vert^{2}\,\mathrm{Im}(\tau)^{k}\,e^{-\frac{4\pi m\mathrm{Im}(z)^{2}}{\mathrm{Im}(\tau)}}
\qquad(\tau\in\mathbb{H},\,z\in\mathbb{C}).
\end {align*}
Let $F$ be a Siegel cusp form of weight $k$ for the Siegel modular group $\mathrm{Sp}_{4}(\mathbb{Z})$, and let $\{\phi_{m}\}_{m\geq 1}$ be the set 
of Jacobi cusp forms appearing in the Fourier--Jacobi expansion of $F$, i.\,e., $\phi_{m}\in J_{k,m}^{\mathrm{cusp}}(\Gamma_{0})$. Then, for any 
$\epsilon>0$, W.~Kohnen proved the following $L^{\infty}$-norm bound for the pointwise Petersson norm of $\phi_{m}$ in~\cite{ko}
\begin{align}
\label{ko:estimate}
\Vert\phi_{m}\Vert_{L^{\infty}}=\sup_{(\tau,z)\in\mathbb{H}\times\mathbb{C}}\Vert\phi_{m}(\tau,z)\Vert_{\mathrm{Pet}}=O_{F,\epsilon}\big(m^{\frac{k}{2}-
\frac{2}{9}+\epsilon}\big),
\end{align}
where the implied constant depends on the Siegel cusp form $F$ and the choice of $\epsilon>0$. Motivated by the Ramanujan--Petersson conjecture,
W.~Kohnen then conjectured the $L^{\infty}$-norm bound
\begin{align*}
\Vert\phi_{m}\Vert_{L^{\infty}}=\sup_{(\tau,z)\in\mathbb{H}\times\mathbb{C}}\Vert\phi_{m}(\tau,z)\Vert_{\mathrm{Pet}}=O_{F,\epsilon}\big(m^{\frac{k-1}
{2}+\epsilon}\big),
\end{align*}
where the implied constant depends on the Siegel cusp form $F$ and the choice of $\epsilon>0$. 

More recently, P.~Anamby and S.~Das established in~\cite{das} a general $L^{\infty}$-norm bound for the pointwise Petersson norm of any Jacobi 
cusp form $\phi\in J_{k,m}^{\mathrm{cusp}}(\Gamma_{0})$, which is $L^{2}$-normalized with respect to the Petersson inner product, i.\,e., for which 
we have
\begin{align*}
\Vert\phi\Vert^{2}_{L^{2}}=\int\limits\limits_{\Gamma_{0}\ltimes\mathbb{Z}^{2}\backslash\mathbb{H}\times\mathbb{C}}\Vert\phi(\tau,z)\Vert^{2}_
{\mathrm{Pet}}\,\frac{\mathrm{d}\xi\wedge\mathrm{d}\eta\wedge\mathrm{d}x\wedge\mathrm{d}y}{\eta^3}=1\qquad(\tau=\xi+i\eta,\,z=x+iy).
\end{align*}
Their $L^{\infty}$-norm bound is (see Theorem~1.4 in\cite{das})
\begin{align}
\label{das:estimate}
\Vert\phi\Vert_{L^{\infty}}=\sup_{(\tau,z)\in\mathbb{H}\times\mathbb{C}}\Vert\phi(\tau,z)\Vert_{\mathrm{Pet}}=O_{\epsilon}\big((k\,m)^{1+\epsilon}
\big),
\end{align}
where the implied constant depends on the choice of $\epsilon>0$.

\subsection{Main results}
\label{subsec-1.2} 
The goal of this article is to provide $L^{\infty}$-norm bounds for the pointwise Petersson norm for Jacobi forms of weight $k\in\mathbb{Z}_{\ge 5}$ 
and index  $m\in\mathbb{Z}_{\ge 1}$ for $\Gamma_{0}$, which are $L^{2}$-normalized. The main result in this respect is given in Theorem~\ref
{thm11} and states for any $L^{2}$-normalized Jacobi cusp form $\phi\in J_{k,m}^{\mathrm{cusp}}(\Gamma_{0})$ the $L^{\infty}$-norm bound
\begin{align}
\label{mainres}
\Vert\phi\Vert_{L^{\infty}}=\sup_{(\tau,z)\in\mathbb{H}\times \mathbb{C}}\Vert\phi(\tau,z)\Vert_{\mathrm{Pet}}=O_{\Gamma_{0},\epsilon}\big(k\,m^
{\frac{7}{4}+\epsilon}\big)
\end{align}
holds, where the implied constant depends on $\Gamma_{0}$ and the choice of $\epsilon>0$. For the proof, we essentially use the representation 
of the Jacobi cusp forms under consideration as linear combinations of modular forms of weight $k-\frac{1}{2}$ for some congruence subgroup of 
$\Gamma_{0}$ (depending on $m$) and suitable Jacobi theta functions; we then need to derive bounds for the $L^{\infty}$-norms of these functions 
to arrive at our result. Comparing our $L^{\infty}$-norm bound~\eqref{mainres} with the $L^{\infty}$-norm bound~\eqref{das:estimate} by P.~Anamby 
and S.~Das, shows that the polynomial growth in $k$ is the same, while the polynomial growth of their bound in $m$ is better by a factor of $m^{\frac
{3}{4}}$. 

In order to be able to derive our $L^{\infty}$-norm bound~\eqref{mainres}, we need $L^{\infty}$-norm bounds for the pointwise Petersson norm 
of cusp forms of positive real weight $k$ and character $\chi$ for any cofinite Fuchsian subgroup $\Gamma$ of $\mathrm{SL}_{2}(\mathbb{R})$. 
Such $L^{\infty}$-norm bounds have been derived  for finite index subgroups of $\mathrm{SL}_{2}(\mathbb{Z})$ by Steiner in~\cite{steiner}. 
Although these $L^{\infty}$-norm bounds are uniform with respect to the index of the finite index subgroups of $\mathrm{SL}_{2}(\mathbb{Z})$, 
there are mild restrictions for a direct application of the $L^{\infty}$-norm bounds derived in~\cite{steiner} to our setting. On the other hand, such 
$L^{\infty}$-norm bounds could also be derived from~\cite{fjk2} with some extra work. However, we provide here new, alternative proofs for the 
results of~\cite{fjk2} applying to any positive real weight $k$ and any character $\chi$ by using the Bergman kernel for the modular curve 
associated to $\Gamma$. 

More specifically, given $\Gamma\subset\mathrm{SL}_{2}(\mathbb{R})$ a Fuchsian subgroup, $k\in\mathbb{R}_{>0}$, and $\chi\colon\Gamma
\rightarrow\mathbb{C}^{\times}$ a character, we let $S_{k,\chi}(\Gamma)$ denote the space of cusp forms of weight $k$ and character $\chi$ for 
$\Gamma$. Denoting by $d_{k}$ the dimension of $S_{k,\chi}(\Gamma)$ and letting $\{f_{1},\ldots,f_{d_{k}}\}$ be an orthonormal basis of $S_{k,
\chi}(\Gamma)$ with respect to the Petersson inner product, the Bergman kernel associated to $S_{k,\chi}(\Gamma)$ is then defined by
\begin{align*}
B_{k,\chi}(\tau,\tau')\coloneqq\sum_{j=1}^{d_{k}}f_{j}(\tau)\overline{f_{j}(\tau')};
\end{align*}
it is straightforward that this definition does not depend on the choice of an orthonormal basis of $S_{k,\chi}(\Gamma)$. The pointwise Petersson 
norm of the Bergman kernel is defined by
\begin{align*}
\Vert B_{k,\chi}(\tau,\tau')\Vert_{\mathrm{Pet}}=\vert B_{k,\chi}(\tau,\tau')\vert\,(\mathrm{Im}(\tau)\mathrm{Im}(\tau'))^{\frac{k}{2}},
\end{align*}
which gives on the diagonal
\begin{align}
\label{aver}
\Vert B_{k,\chi}(\tau,\tau)\Vert_{\mathrm{Pet}}=\sum\limits_{j=1}^{d_{k}}\Vert f_{j}(\tau)\Vert^{2}_{\mathrm{Pet}}.
\end{align}
As a second main result of this article, we establish in Theorem~\ref{thm4}, assuming that $k\in\mathbb{R}_{\ge 5}$, for $\Gamma$ being 
cocompact without elliptic elements the $L^{\infty}$-norm bound
\begin{align*}
\Vert B_{k,\chi}\Vert_{L^{\infty}}=\sup_{\tau\in\mathbb{H}}\Vert B_{k,\chi}(\tau,\tau)\Vert_{\mathrm{Pet}}=O_{\Gamma}(k);
\end{align*}
moreover, for $\Gamma$ being cofinite, we give the $L^{\infty}$-norm bound
\begin{align*}
\Vert B_{k,\chi}\Vert_{L^{\infty}}=\sup_{z\in\mathbb{H}}\Vert B_{k,\chi}(\tau,\tau)\Vert_{\mathrm{Pet}}=O_{\Gamma}\big(k^{\frac{3}{2}}\big),
\end{align*}
where the implied constants depend only on the Fuchsian subgroup~$\Gamma$. Due to the relation~\eqref{aver}, these results reprove the 
${L^{\infty}}$-norm bounds on average obtained in~\cite{fjk2}, but now for any real weight $k\in\mathbb{R}_{\ge 5}$ and any character $\chi$. 
Based on these results, we are then also able to prove the uniformity of the above ${L^{\infty}}$-norm bounds with respect to the subgroup 
$\Gamma$ in Theorem~\ref{thm6}.

\subsection{Outline}
\label{subsec-1.3} 
Let us briefly outline the contents of this article. In the subsequent, second section we collect all the necessary prerequisites for the sequel of 
the paper. In particular, we introduce the definitions of cusp forms and Jacobi cusp forms together with their respective pointwise Petersson
norms. Furthermore, we define the Bergman kernel for modular curves and state its basic properties. 

The third section is devoted to the revisiting of the ${L^{\infty}}$-norm bounds on average obtained in~\cite{fjk2}, but now for any real weight 
$k\in\mathbb{R}_{\ge 5}$ and any character $\chi$. Here, the proofs of Theorem~\ref{thm4} and Theorem~\ref{thm6} are provided.

In the fourth section, the ${L^{\infty}}$-norm bound~\eqref{mainres} is proven in Theorem~\ref{thm11}. In order to establish this result, we 
combine the results of the third section applied to modular forms of half-integral weight with $L^{\infty}$-norm bounds of classical Jacobi 
theta functions, which are derived by directly estimating their defining series expansions.

\section{Preliminaries}
\label{sec-2}

\subsection{Hyperbolic metric}
\label{subsec-2.1}
Let $\mathbb{H}\coloneqq\{\tau\in\mathbb{C}\,\vert\,\tau=\xi+i\eta,\,\eta>0\}$ be the upper half-plane. We denote by $\mathrm{d}s^{2}_{\mathrm
{hyp}}(\tau)$ the line element and by $\mu_{\mathrm{hyp}}(\tau)$ the volume form corresponding to the hyperbolic metric on $\mathbb{H}$, which 
is compatible with the complex structure of $\mathbb{H}$ and has constant curvature equal to $-1$. Locally on $\mathbb{H}$, we have
\begin{align*}
\mathrm{d}s^{2}_{\mathrm{hyp}}(\tau)=\frac{\mathrm{d}\xi^{2}+\mathrm{d}\eta^{2}}{\eta^{2}}\quad\textrm{and}\quad\mu_{\mathrm{hyp}}(\tau)=
\frac{\mathrm{d}\xi\wedge\mathrm{d}\eta}{\eta^{2}}\,.
\end{align*}
For $\tau,\tau'\in\mathbb{H}$, we let $\mathrm{dist}_{\mathrm{hyp}}(\tau,\tau')$ denote the hyperbolic distance between these two points. For
later purposes, it is useful to introduce the displacement function
\begin{align}
\label{dist-eqn}
\sigma(\tau,\tau')\coloneqq\cosh^{2}\bigg(\frac{\mathrm{dist}_{\mathrm{hyp}}(\tau,\tau')}{2}\bigg)=\frac{\vert\tau-\overline{\tau}'\vert^{2}}{4\,\mathrm
{Im}(\tau)\mathrm{Im}(\tau')}\,.
\end{align}

\subsection{Quotient space}
\label{subsec-2.2}
Let $\Gamma\subset\mathrm{SL}_{2}(\mathbb{R})$ be a Fuchsian subgroup acting by fractional linear transformations on $\mathbb{H}$. Let $X_
{\Gamma}$ be the quotient space $\Gamma\backslash\mathbb{H}$ and $g_{\Gamma}$ the genus of $X_{\Gamma}$. In the sequel, we identify 
$X_{\Gamma}$ with a fundamental domain $\mathcal{F}_{\Gamma}\subset\mathbb{H}$ for the group $\Gamma$, which we assume to be closed 
and connected.

Denote by 
\begin{align*}
\mathcal{P}_{\Gamma}=\{p_{1},\ldots,p_{s}\}
\end{align*}
the set of cusps of $\mathcal{F}_{\Gamma}$. Let $\sigma_{\mathcal{P},j}\in\mathrm{SL}_{2}(\mathbb{R})$ be a scaling matrix of the cusp $p_{j}$, 
that is, $p_{j}=\sigma_{\mathcal{P},j}i\infty$ with stabilizer subgroup $\Gamma_{p_{j}}$ described as
\begin{align}
\label{parascaling}
\sigma_{\mathcal{P},j}^{-1}\Gamma_{p_{j}}\sigma_{\mathcal{P},j}=
\begin{cases}
\bigg\langle\begin{pmatrix}1&1\\0&1\end{pmatrix}\bigg\rangle,&\text{ if }-\mathrm{id}\not\in\Gamma, \\
\bigg\langle\pm\begin{pmatrix}1&1\\0&1\end{pmatrix}\bigg\rangle,&\text{ if }-\mathrm{id}\in\Gamma,
\end{cases}
\qquad(j=1,\ldots,s).
\end{align}
 For $Y>0$, we let $\mathcal{F}^{Y}_{j}\subset\mathcal{F}_{\Gamma}$ denote the neighborhood of the cusp $p_{j}$ characterized by
\begin{align*}
\sigma_{\mathcal{P},j}^{-1}\mathcal{F}^{Y}_{j}=\{\tau=\xi+i\eta\in\mathbb{H}\,\vert\,-1/2\le\xi\le 1/2,\,\eta\ge Y\}\qquad(j=1,\ldots,s).
\end{align*}
With these notations, we define $\mathcal{F}_{Y}$ to be the closure of the complement of the union $\mathcal{F}^{Y}_{1}\cup\ldots\cup\mathcal
{F}^{Y}_{s}$ in $\mathcal{F}_{\Gamma}$, i.\,e., 
\begin{align}
\label{calfy}
\mathcal{F}_{Y}\coloneqq\mathrm{cl}\big(\mathcal{F}_{\Gamma}\setminus\big(\mathcal{F}^{Y}_{1}\cup\ldots\cup\mathcal{F}^{Y}_{s}\big)\big),
\end{align}
which is compact. For $Y$ being sufficiently large, we note that $\mathcal{F}_{Y}=\mathcal{F}_{\Gamma}$, if $\Gamma$ is cocompact. We 
choose $0<m_{Y}<M_{Y}$ such that for all $\tau\in\mathcal{F}_{Y}$ the inequalities
\begin{align*}
m_{Y}\le\mathrm{Im}(\sigma_{\mathcal{P},j}^{-1}\tau)\le M_{Y}
\end{align*}
hold for all $j=1,\ldots,s$; we note that $m_{Y}$ and $M_{Y}$ depend on the choice of $Y$.

Denote by 
\begin{align*}
\mathcal{E}_{\Gamma}=\{e_{1},\ldots,e_{t}\}
\end{align*}
the set of elliptic fixed points of $\mathcal{F}_{\Gamma}$. Let $\Gamma_{e_{j}}$ and $m_{j}$ denote the stabilizer subgroup and order of the 
elliptic fixed point $e_{j}$, respectively.

We denote the hyperbolic length of the shortest closed geodesic on $X_{\Gamma}$ by $\ell_{\Gamma}$. For a domain $D\subset\mathbb{H}$, we 
denote its hyperbolic diameter by $\mathrm{diam}_{\mathrm{hyp}}(D)$ and its hyperbolic volume by $\mathrm{vol}_{\mathrm{hyp}}(D)$. Finally, the
injectivity radius $r_{\Gamma}$ is defined by
\begin{align}
\label{irad}
r_{\Gamma}\coloneqq\inf\bigg\{\mathrm{dist}_{\mathrm{hyp}}(\tau,\gamma\tau)\,\bigg\vert\,\tau\in\mathcal{F}_{\Gamma},\,\gamma\in\Gamma
\setminus\bigg(\bigcup_{j=1}^{s}\Gamma_{p_{j}}\cup\bigcup_{j=1}^{t}\Gamma_{e_{j}}\cup\{\mathrm{id}\}\bigg)\bigg\}.
\end{align}
We note that if $X_{\Gamma}$ is compact without elliptic fixed points, i.\,e., $\mathcal{P}_{\Gamma}=\mathcal{E}_{\Gamma}=\emptyset$, then the
injectivity radius $r_{\Gamma}$ equals the length of the shortest closed geodesic $\ell_{\Gamma}$ of $X_{\Gamma}$.

\subsection{Cusp forms and Bergman kernel}
\label{subsec-2.3}
For $k\in\mathbb{R}_{>0}$ and a character $\chi\colon\Gamma\rightarrow\mathbb{C}^{\times}$, we let $S_{k,\chi}(\Gamma)$ denote the space of 
cusp forms of weight $k$ and character $\chi$ for $\Gamma$, i.\,e., the space of holomorphic functions $f\colon\mathbb{H}\rightarrow\mathbb{C}$, 
which have the transformation behavior
\begin{align*}
f(\gamma z)(cz+d)^{-k}=\chi(\gamma)f(z)
\end{align*}
for all $\gamma=\big(\begin{smallmatrix}a&b\\c&d\end{smallmatrix}\big)\in\Gamma$, and which vanish at all the cusps of $\mathcal{F}_{\Gamma}$. 
Given $f\in S_{k,\chi}(\Gamma)$, we define
\begin{align*}
\Vert f(\tau)\Vert_{\mathrm{Pet}}^{2}\coloneqq\vert f(\tau)\vert^{2}\,\eta^{k}\qquad(\tau=\xi+i\eta),
\end{align*}
which defines a $\Gamma$-invariant function on $\mathbb{H}$ called the pointwise Petersson norm of $f$. 

The space $S_{k,\chi}(\Gamma)$ is equipped with the Petersson inner product
\begin{align}
\label{pet-ip-gamma}
\langle f_{1},f_{2}\rangle_{\mathrm{Pet}}\coloneqq\int\limits_{\mathcal{F}_{\Gamma}}f_{1}(\tau)\overline{f_{2}(\tau)}\,\eta^{k}\mu_{\mathrm{hyp}}(\tau)
\qquad(f_{1},f_{2}\in S_{k,\chi}(\Gamma)).
\end{align}
Let $d_{k}$ denote the dimension of $S_{k,\chi}(\Gamma)$ and let $\{f_{1},\ldots,f_{d_{k}}\}$ be an orthonormal basis of $S_{k,\chi}(\Gamma)$ with 
respect to the Petersson inner product. Then, the Bergman kernel associated to $S_{k,\chi}(\Gamma)$ is defined by
\begin{align*}
B_{k,\chi}(\tau,\tau')\coloneqq\sum_{j=1}^{d_{k}}f_{j}(\tau)\overline{f_{j}(\tau')}.
\end{align*}
It is obvious that this definition does not depend on the choice of an orthonormal basis of $S_{k,\chi}(\Gamma)$.

The Bergman kernel $B_{k,\chi}(\tau,\tau')$ is a holomorphic cusp form of weight $k$ and character $\chi$ for $\Gamma$ in the $\tau$-variable, and 
an anti-holomorphic cusp form of weight $k$ and character $\overline{\chi}$ for $\Gamma$ in the $\tau'$-variable. Hence, the pointwise Petersson 
norm of the Bergman kernel is given by
\begin{align*}
\Vert B_{k,\chi}(\tau,\tau')\Vert_{\mathrm{Pet}}=\vert B_{k,\chi}(\tau,\tau')\vert\,(\eta\eta')^{\frac{k}{2}},
\end{align*}
which is a $\Gamma$-invariant function on $\mathbb{H}\times\mathbb{H}$ with respect to both variables.  

Moreover, $B_{k,\chi}(\tau,\tau')$ is the reproducing kernel for $S_{k,\chi}(\Gamma)$, i.\,e., we have
\begin{align*}
\int\limits_{\mathcal{F}_{\Gamma}}B_{k,\chi}(\tau,\tau')f(\tau')\eta'^{k}\mu_{\mathrm{hyp}}(\tau')=f(\tau)
\qquad(\tau'=\xi'+i\eta')
\end{align*}
for any $f\in S_{k,\chi}(\Gamma)$. Therefore, for $k\in\mathbb{R}_{>3}$, the Bergman kernel $B_{k,\chi}(\tau,\tau')$ can also be represented in the 
following form (see Proposition~1.3 on p.~77 in~\cite{freitag})
\begin{align}
\label{bkseries}
B_{k,\chi}(\tau,\tau')=\frac{(2i)^{k}(k-1)}{4\pi}\sum_{\gamma=\left(\begin{smallmatrix}a&b\\c&d\end{smallmatrix}\right)\in\Gamma}\frac{1}{(\tau-\gamma
\overline{\tau}')^{k}}\frac{1}{\chi(\gamma)(c\overline{\tau}'+d)^{k}}\,.
\end{align}
Note that the formula for the Bergman kernel given in~\cite{freitag} is missing a factor of $(2i)^{k}$.

\subsection{Counting function}
\label{subsec-2.4}
Given $\tau\in\mathbb{H}$ and $\rho\in\mathbb{R}_{\ge 0}$, we recall from~\cite{jl} the counting function
\begin{align*}
N_{\Gamma}(\tau;\rho)\coloneqq\vert\mathcal{N}_{\Gamma}(\tau;\rho)\vert,
\end{align*}
where
\begin{align*}
\mathcal{N}_{\Gamma}(\tau;\rho)\coloneqq\bigg\{\gamma\in\Gamma\setminus\bigg(\bigcup_{j=1}^{s}\Gamma_{p_{j}}\cup\bigcup_{j=1}^{t}\Gamma_
{e_{j}}\cup\{\mathrm{id}\}\bigg)\,\bigg\vert\,\mathrm{dist}_{\mathrm{hyp}}(\tau,\gamma\tau)\le\rho\bigg\}.
\end{align*} 
Let now $f$ be a positive, smooth, and decreasing function on $\mathbb{R}_{\ge 0}$. Then, adapting the arguments from~\cite{jl} to Fuchsian 
subgroups of $\mathrm{SL}_{2}(\mathbb{R})$, we have for any  $\tau\in\mathbb{H}$ and any $\delta\ge r_{\Gamma}/2$ the inequality 
\begin{align}
\label{jlineq1}
\int\limits_{0}^{\infty}f(\rho)\,\mathrm{d}N_{\Gamma}(\tau;\rho)\le&\int\limits_{0}^{\delta}f(\rho)\,\mathrm{d}N_{\Gamma}(\tau;\rho)+\frac{2\vert
\mathrm{Cent}(\Gamma)\vert\cosh(r_{\Gamma}/4)}{\sinh(r_{\Gamma}/4)}\sinh(\delta)f(\delta) \\[2mm]
\notag
&+\frac{\vert\mathrm{Cent}(\Gamma)\vert}{2\sinh^{2}(r_{\Gamma}/4)}\int\limits_{\delta}^{\infty}f(\rho)\sinh(\rho+r_{\Gamma}/2\big)\,\mathrm{d}
\rho;
\end{align}  
here $\mathrm{Cent}(\Gamma)$ denotes the center of $\Gamma$. Note that our definition~\eqref{irad} of injectivity radius differs from the one 
used in~\cite{jl} by a factor of $2$, and the inequality~\eqref{jlineq1} takes this fact into account. 

\subsection{Jacobi forms}
\label{subsec-2.5}
For $k, m\in\mathbb{N}$, we let $J_{k,m}^{\mathrm{cusp}}(\Gamma_{0})$ denote the space of Jacobi cusp forms of weight $k$ and index $m$ for 
$\Gamma_{0}=\mathrm{SL}_{2}(\mathbb{Z})$, i.\,e., the space of holomorphic functions $\phi\colon\mathbb{H}\times\mathbb{C}\rightarrow\mathbb
{C}$, which have the transformation behaviour 
\begin{align*}
\phi\Big(\frac{a\tau+b}{c\tau+d},\frac{z+\lambda\tau+\mu}{c\tau+d}\Big)(c\tau+d)^{-k}\exp\Big(2\pi im\Big(\lambda^{2}\tau+2\lambda z-\frac{c(z+
\lambda\tau+\mu)^{2}}{c\tau+d}\Big)\Big)=\phi(\tau,z)
\end{align*}
for all $\big[\big(\begin{smallmatrix}a&b\\c&d\end{smallmatrix}\big),(\lambda,\mu)\big]\in\Gamma_{0}\ltimes\mathbb{Z}^{2}$, and which have a 
Fourier expansion of the form
\begin{align*}
\phi(\tau,z)=\sum\limits_{\substack{n\in\mathbb{N},\,r\in\mathbb{Z}\\4mn-r^{2}>0}}c(n,r)\,q^{n}\zeta^{r}\qquad(q=e^{2\pi i\tau},\,\zeta=e^{2\pi iz}).
\end{align*}
Given $\phi\in J_{k,m}^{\mathrm{cusp}}(\Gamma_{0})$, we define
\begin{align*}
\Vert\phi(\tau,z)\Vert_{\mathrm{Pet}}^{2}\coloneqq\vert\phi(\tau,z)\vert^{2}\,\eta^{k}\,e^{-\frac{4\pi my^{2}}{\eta}}
\qquad(\tau=\xi+i\eta,\,z=x+iy),
\end{align*}
which defines a $\Gamma_{0}\ltimes\mathbb{Z}^{2}$-invariant function on $\mathbb{H}\times\mathbb{C}$ called the pointwise Petersson norm 
of~$\phi$. 

Let $\mathcal{D}_{\Gamma_{0}}$ denote a fundamental domain of the quotient space $Y_{\Gamma_{0}}\coloneqq\Gamma_{0}\ltimes\mathbb{Z}^
{2}\backslash\mathbb{H}\times\mathbb{C}$, which is a $2$-dimensional complex manifold. The space $ J_{k,m}^{\mathrm{cusp}}(\Gamma_{0})$ is 
equipped with the Petersson inner product
\begin{align}
\label{pet-ip}
\langle\phi_{1},\phi_{2}\rangle_{\mathrm{Pet}}\coloneqq\int\limits\limits_{\mathcal{D}_{\Gamma_{0}}}\phi_{1}(\tau,z)\overline{\phi_{2}(\tau,z)}\,\eta^
{k}\,e^{-\frac{4\pi my^{2}}{\eta}}\,\frac{\mathrm{d}\xi\wedge\mathrm{d}\eta\wedge\mathrm{d}x\wedge\mathrm{d}y}{\eta^3}\qquad(\phi_{1},\phi_{2}\in 
J_{k,m}^{\mathrm{cusp}}(\Gamma_{0})).
\end{align}

For $\phi\in J_{k,m}^{\mathrm{cusp}}(\Gamma_{0})$, one has the decomposition
\begin{align}
\label{eich-zag1}
\phi(\tau,z)=\sum_{\mu=0}^{2m-1}h_{\mu}(\tau)\vartheta_{\mu,m}(\tau,z),
\end{align}
where the function $h_{\mu}$ by Remark~2 on p.~287 of~\cite{k2} is a cusp form of weight $(k-\frac{1}{2})$ for the finite index subgroup 
\begin{align*}
\Gamma_{1}\coloneqq\Gamma_{0,1}(4m)=\bigg\{\begin{pmatrix}a&b\\c&d\end{pmatrix}\in\Gamma_{0}\,\bigg\vert\,c\equiv 0,\,d\equiv 1\mod 4m
\bigg\}\subseteq\Gamma_{0}(4m),
\end{align*}
and $\vartheta_{\mu,m}$ is the Jacobi theta function
\begin{align}
\label{defn:theta}
\vartheta_{\mu,m}(\tau,z)\coloneqq\sum_{n\in\mathbb{Z}}e^{2\pi im\tau\big(n-\frac{\mu}{2m}\big)^{2}+2\pi iz(2mn-\mu)}.
\end{align}
In fact, it is shown in Theorem~5.1 of~\cite{eich-zag} that the decomposition~\eqref{defn:theta} gives rise to the isomorphism
\begin{align*}
J_{k,m}^{\mathrm{cusp}}(\Gamma_{0})\cong\mathcal{V}_{k-\frac{1}{2}}(\Gamma_{0}),
\end{align*} 
where $\mathcal{V}_{k-\frac{1}{2}}(\Gamma_{0})$ denotes the complex vector space of vector-valued cusp forms of weight $(k-\frac{1}{2})$ with 
suitable transformation behaviour with respect to $\Gamma_{0}$.

Let now
\begin{align*}
\phi_{1}(\tau,z)=\sum_{\mu=0}^{2m-1}h_{\mu,1}(\tau)\vartheta_{\mu,m}(\tau,z)\quad\text{and}\quad\phi_{2}(\tau,z)=\sum_{\mu=0}^{2m-1}h_{\mu,2}
(\tau)\vartheta_{\mu,m}(\tau,z)
\end{align*}
be two Jacobi cusp forms of weight $k$ and index $m$ for $\Gamma_{0}$. Then, the decomposition~\eqref{eich-zag1} gives rise to the equality 
(see Theorem~5.3 in~\cite{eich-zag})
\begin{align}
\label{pet-ip2}
\langle\phi_{1},\phi_{2}\rangle_{\mathrm{Pet}}=\frac{1}{\sqrt{4m}}\int\limits_{\mathcal{F}_{\Gamma_{0}}}\sum_{\mu=0}^{2m-1}h_{\mu,1}(\tau)
\overline{h_{\mu,2}(\tau)}\,\eta^{k-\frac{1}{2}}\,\frac{\mathrm{d}\xi\wedge\mathrm{d}\eta}{\eta^2},
\end{align}
where $\mathcal{F}_{\Gamma_{0}}$ denotes a fundamental domain for the quotient space $X_{\Gamma_{0}}=\Gamma_{0}\backslash\mathbb{H}$.

\section{$L^{\infty}$-norm bounds for cusp forms revisited}
\label{sec3}

Refining arguments of~\cite{am1} and~\cite{am2}, we first derive bounds for the Bergman kernel along the diagonal. 
\begin{prop}
\label{prop1}
With notations as above, let $\Gamma\subset\mathrm{SL}_{2}(\mathbb{R})$ be a cocompact Fuchsian subgroup without elliptic elements. Then,
for $k\in\mathbb{R}_{\ge 5}$ and $\tau\in\mathbb{H}$, we have the bound
\begin{align*}
\Vert B_{k,\chi}(\tau,\tau)\Vert_{\mathrm{Pet}}\le\frac{k-1}{2\pi}+\frac{3(k-1)}{\pi\cosh^{k-4}(r_{\Gamma}/4)}\bigg(1+\frac{1}{\sinh^{2}(r_{\Gamma}/4)}
\bigg).
\end{align*}
\begin{proof}
Letting $k\in\mathbb{R}_{\ge 5}$ and considering the Bergman kernel~\eqref{bkseries} on the diagonal, we derive by means of relation~\eqref
{dist-eqn} the bound
\begin{align}
\notag
&\Vert B_{k,\chi}(\tau,\tau)\Vert_{\mathrm{Pet}}=\frac{2^{k}(k-1)}{4\pi}\bigg\vert\sum_{\gamma=\left(\begin{smallmatrix}a&b\\c&d\end{smallmatrix}
\right)\in\Gamma}\frac{1}{(\tau-\gamma\overline{\tau})^{k}}\frac{1}{\chi(\gamma)(c\overline{\tau}+d)^{k}}\bigg\vert\,\mathrm{Im}(\tau)^{k} \\
\notag
&\qquad\le\frac{k-1}{4\pi}\sum_{\gamma\in\Gamma}\bigg(\frac{4\,\mathrm{Im}(\tau)\mathrm{Im}(\gamma\tau)}{\vert\tau-\gamma\overline{\tau}
\vert^{2}}\bigg)^{k/2}=\frac{k-1}{4\pi}\sum_{\gamma\in\Gamma}\frac{1}{\cosh^{k}(\mathrm{dist}_{\mathrm{hyp}}(\tau,\gamma\tau)/2)} \\[1mm]
\label{prop1-eqn1}
&\qquad=\frac{k-1}{4\pi}\bigg(\vert\mathrm{Cent}(\Gamma)\vert+\sum_{\gamma\in\Gamma\setminus\mathrm{Cent}(\Gamma)}\frac{1}{\cosh^{k}
(\mathrm{dist}_{\mathrm{hyp}}(\tau,\gamma\tau)/2)}\bigg).
\end{align}
Substituting $\delta=r_{\Gamma}/2$ in inequality~\eqref{jlineq1} and using the fact that $\vert\mathrm{Cent}(\Gamma)\vert\le 2$, we derive
\begin{align}
\notag
&\sum_{\gamma\in\Gamma\setminus\mathrm{Cent}(\Gamma)}\frac{1}{\cosh^{k}(\mathrm{dist}_{\mathrm{hyp}}(\tau,\gamma\tau)/2)} \\
\label{prop1-eqn2}
&\qquad\le\int\limits_{0}^{r_{\Gamma}/2}\frac{\mathrm{d}N_{\Gamma}(\tau;\rho)}{\cosh^{k}(\rho/2)}+\frac{8}{\cosh^{k-2}(r_{\Gamma}/4)}+\frac{1}
{\sinh^{2}(r_{\Gamma}/4)}\int\limits_{r_{\Gamma}/2}^{\infty}\frac{\sinh(\rho+r_{\Gamma}/2)}{\cosh^{k}(\rho/2)}\,\mathrm{d}\rho.
\end{align}
From the defining equation~\eqref{irad} of the injectivity radius $r_{\Gamma}$, we find for the first term of~\eqref{prop1-eqn2} that
\begin{align}
\label{prop1-eqn3}
\int\limits_{0}^{r_{\Gamma}/2}\frac{\mathrm{d}N_{\Gamma}(\tau;\rho)}{\cosh^{k}(\rho/2)}=0.
\end{align} 
With regard to the third term of~\eqref{prop1-eqn2}, we recall the bound (12) in~\cite{am1}, which states for any $k\in\mathbb{R}_{\ge 5}$ and 
any $\delta\ge 0$ (note that we have replaced $2k$ by $k$) that
\begin{align}
\notag
&\frac{1}{\sinh^{2}(r_{\Gamma}/4)}\int\limits_{\delta}^{\infty}\frac{\sinh(\rho+r_{\Gamma}/2)}{\cosh^{k}(\rho/2)}\,\mathrm{d}\rho \\
\label{prop1-eqn4}
&\quad\le\frac{4}{(k-2)\cosh^{k-2}(\delta/2)}\bigg(2+\frac{1}{\sinh^{2}(r_{\Gamma}/4)}\bigg)+\frac{8}{(k-4)\cosh^{k-4}(\delta/2)}\cdot\frac{1}{\sinh^
{2}(r_{\Gamma}/4)}.
\end{align}
From the elementary inequality $\cosh^{k-4}(r_{\Gamma}/4)\le\cosh^{k-2}(r_{\Gamma}/4)$ and recalling that $k\in\mathbb{R}_{\ge 5}$, we now 
derive from~\eqref{prop1-eqn4} with $\delta=r_{\Gamma}/2$ the bound
\begin{align}
\notag
&\frac{1}{\sinh^{2}(r_{\Gamma}/4)}\int\limits_{r_{\Gamma}/2}^{\infty}\frac{\sinh(\rho+r_{\Gamma}/2)}{\cosh^{k}(\rho/2)}\,\mathrm{d}\rho \\
\notag
&\le\frac{4}{(k-2)\cosh^{k-2}(r_{\Gamma}/4)}\bigg(2+\frac{1}{\sinh^{2}(r_{\Gamma}/4)}\bigg)+\frac{8}{(k-4)\cosh^{k-4}(r_{\Gamma}/4)}\cdot\frac{1}
{\sinh^{2}(r_{\Gamma}/4)} \\[2mm]
\label{prop1-eqn5}
&\le\frac{4}{\cosh^{k-4}(r_{\Gamma}/4)}\bigg(1+\frac{1}{\sinh^{2}(r_{\Gamma}/4)}\bigg)+\frac{8}{(k-4)\cosh^{k-4}(r_{\Gamma}/4)}\cdot\frac{1}
{\sinh^{2}(r_{\Gamma}/4)}.
\end{align}
Combining the bounds~\eqref{prop1-eqn1},~\eqref{prop1-eqn2} with~\eqref{prop1-eqn3},~\eqref{prop1-eqn5}, and using the fact that $k\in
\mathbb{R}_{\ge 5}$, completes the proof of the proposition.  
\end{proof}
\end{prop}

\begin{prop}
\label{prop3}
With notations as above, let $\Gamma\subset\mathrm{SL_{2}(\mathbb{R})}$ be a cofinite Fuchsian subgroup. Then, for $k\in\mathbb{R}_{\ge5}$ 
and $\tau\in\mathbb{H}$, we have the bound
\begin{align*}
\Vert B_{k,\chi}(\tau,\tau)\Vert_{\mathrm{Pet}}&\le\frac{k-1}{2\pi}+\frac{3(k-1)}{\pi\cosh^{k-4}(r_{\Gamma}/4)}\bigg(1+\frac{1}{\sinh^{2}(r_{\Gamma}/
4)}\bigg) \\[2mm]
&\quad+\frac{k-1}{4\pi}\sum_{e_{j}\in\mathcal{E}_{\Gamma}}(m_{j}-1)+\frac{2(k-1)}{\sqrt{\pi}}\cdot\frac{\Gamma((k-1)/2)}{\Gamma(k/2)}\sum_
{p_{j}\in\mathcal{P}_{\Gamma}}\mathrm{Im}(\sigma_{\mathcal{P},j}^{-1}\tau),
\end{align*}
where $\sigma_{\mathcal{P},j}$ is the scaling matrix associated to the cusp $p_{j}\in\mathcal{P}_{\Gamma}$ defined in~\eqref{parascaling}.
\begin{proof}
For $k\in\mathbb{R}_{\ge 5}$ and $\tau\in\mathbb{H}$, using the bound~\eqref{prop1-eqn1} and the fact that $\vert\mathrm{Cent}(\Gamma)\vert
\le 2$, we derive
\begin{align}
\notag
\Vert B_{k,\chi}(\tau,\tau)\Vert_{\mathrm{Pet}}&\le\frac{k-1}{2\pi}+\frac{k-1}{4\pi}\sum_{\gamma\in\Gamma\setminus\big(\bigcup_{e_{j}\in\mathcal
{E}_{\Gamma}}\Gamma_{e_{j}}\,\cup\,\bigcup_{p_{j}\in\mathcal{P}_{\Gamma}}\Gamma_{p_{j}}\big)}\frac{1}{\cosh^{k}(\mathrm{dist}_{\mathrm{hyp}}
(\tau,\gamma\tau)/2)} \\[1mm]
\notag
&\quad+\frac{k-1}{4\pi}\sum_{e_{j}\in\mathcal{E}_{\Gamma}}\,\sum_{\gamma\in\Gamma_{e_{j}}\setminus\mathrm{Cent}(\Gamma)}\frac{1}{\cosh^{k}
(\mathrm{dist}_{\mathrm{hyp}}(\tau,\gamma\tau)/2)} \\[1mm]
\label{prop3-eqn1}
&\quad+\frac{k-1}{4\pi}\sum_{p_{j}\in\mathcal{P}_{\Gamma}}\,\sum_{\gamma\in\Gamma_{p_{j}}\setminus\mathrm{Cent}(\Gamma)}\frac{1}{\cosh^{k}
(\mathrm{dist}_{\mathrm{hyp}}(\tau,\gamma\tau)/2)}.
\end{align}
Adapting our arguments from Proposition~\ref{prop1} to the second summand on the right-hand side of~\eqref{prop3-eqn1}, we arrive at the bound
\begin{align}
\notag
\frac{k-1}{4\pi}\sum_{\gamma\in\Gamma\setminus\big(\bigcup_{e_{j}\in\mathcal{E}_{\Gamma}}\Gamma_{e_{j}}\,\cup\,\bigcup_{p_{j}\in\mathcal{P}_
{\Gamma}}\Gamma_{p_{j}}\big)}\frac{1}{\cosh^{k}(\mathrm{dist}_{\mathrm{hyp}}(\tau,\gamma\tau)/2)} \\[2mm]
\label{prop3-eqn2}
\le\frac{3(k-1)}{\pi\cosh^{k-4}(r_{\Gamma}/4)}\bigg(1+\frac{1}{\sinh^{2}(r_{\Gamma}/4)}\bigg).
\end{align}
For the third term on the right-hand side of~\eqref{prop3-eqn1}, we trivially have the bound
\begin{align}
\label{prop3-eqn3}
\frac{k-1}{4\pi}\sum_{e_{j}\in\mathcal{E}_{\Gamma}}\,\sum_{\gamma\in\Gamma_{e_{j}}\setminus\mathrm{Cent}(\Gamma)}\frac{1}{\cosh^{k}(\mathrm
{dist}_{\mathrm{hyp}}(\tau,\gamma\tau)/2)}\le\frac{k-1}{4\pi}\sum_{e_{j}\in\mathcal{E}_{\Gamma}}(m_{j}-1).
\end{align}
From the definition of the scaling matrix~\eqref{parascaling} and using the fact that $\vert\mathrm{Cent(\Gamma)}\vert\le 2$, we find 
\begin{align}
\notag
&\frac{k-1}{4\pi}\sum_{p_{j}\in\mathcal{P}_{\Gamma}}\,\sum_{\gamma\in\Gamma_{p_{j}}\setminus\mathrm{Cent}(\Gamma)}\frac{1}{\cosh^{k}(\mathrm
{dist}_{\mathrm{hyp}}(\tau,\gamma\tau)/2)} \\[1mm]
\label{prop3-eqn4}
&\qquad\le\frac{k-1}{2\pi}\sum_{p_{j}\in\mathcal{P}_{\Gamma}}\,\sum_{n\in\mathbb{Z}\backslash\{0\}}\frac{1}{\cosh^{k}(\mathrm{dist}_{\mathrm{hyp}}
(\sigma_{\mathcal{P},j}^{-1}\tau,\sigma_{\mathcal{P},j}^{-1}\tau+n)/2)}.
\end{align}
We now recall the bound~(18) in~\cite{am1}, which gives for $k\in\mathbb{R}_{\ge 5}$, $p_{j}\in\mathcal{P}_{\Gamma}$, and $\tau,\tau'\in\mathbb{H}$
(note that we have replaced $2k$ by $k$) the bound
\begin{align}
\notag
&\frac{k-1}{2\pi}\sum_{p_{j}\in\mathcal{P}_{\Gamma}}\,\sum_{n\in\mathbb{Z}\backslash\{0\}}\frac{1}{\cosh^{k}(\mathrm{dist}_{\mathrm{hyp}}(\sigma_
{\mathcal{P},j}^{-1}\tau,\sigma_{\mathcal{P},j}^{-1}\tau'+n)/2)} \\[1mm]
\label{prop3-eqn5}
&\qquad\le
\frac{k-1}{\sqrt{\pi}}\cdot\frac{\Gamma\big((k-1)/2\big)}{\Gamma(k/2)}\sum_{p_{j}\in\mathcal{P}_{\Gamma}}\frac{\big(4\,\mathrm{Im}(\sigma_{\mathcal
{P},j}^{-1}\tau)\mathrm{Im}(\sigma_{\mathcal{P},j}^{-1}\tau')\big)^{k/2}}{\big(\mathrm{Im}(\sigma_{\mathcal{P},j}^{-1}\tau)+\mathrm{Im}(\sigma_{\mathcal
{P},j}^{-1}\tau')\big)^{k-1}}.
\end{align}
Substituting $\tau=\tau'$ in~\eqref{prop3-eqn5} and combining it with~\eqref{prop3-eqn4}, we arrive for the fourth term on the right-hand side of~\eqref
{prop3-eqn1} at the bound
\begin{align}
\notag
&\frac{k-1}{4\pi}\sum_{p_{j}\in\mathcal{P}_{\Gamma}}\,\sum_{\gamma\in\Gamma_{p_{j}}\setminus\mathrm{Cent}(\Gamma)}\frac{1}{\cosh^{k}(\mathrm
{dist}_{\mathrm{hyp}}(\tau,\gamma\tau)/2)} \\[1mm]
\label{prop3-eqn6}
&\qquad\le\frac{2(k-1)}{\sqrt{\pi}}\cdot\frac{\Gamma\big((k-1)/2\big)}{\Gamma(k/2)}\sum_{p_{j}\in\mathcal{P}_{\Gamma}}\mathrm{Im}(\sigma_{\mathcal
{P},j}^{-1}\tau).
\end{align}
Combining the bounds~\eqref{prop3-eqn2},~\eqref{prop3-eqn3}, and~\eqref{prop3-eqn6} with~\eqref{prop3-eqn1} completes the proof of the proposition. 
\end{proof}
\end{prop}

\begin{thm}
\label{thm4}
With notations as above, let $\Gamma\subset\mathrm{SL_{2}(\mathbb{R})}$ be a cofinite Fuchsian subgroup and $k\in\mathbb{R}_{\ge 5}$. Then, 
if $\Gamma$ is cocompact without elliptic elements, we have the $L^{\infty}$-norm bound
\begin{align}
\label{thm4:eqn1}
\Vert B_{k,\chi}\Vert_{L^{\infty}}=\sup_{\tau\in\mathbb{H}}\Vert B_{k,\chi}(\tau,\tau)\Vert_{\mathrm{Pet}}=O_{\Gamma}(k).
\end{align}
Moreover, if $\Gamma$ is cofinite, we have the $L^{\infty}$-norm bound
\begin{align}
\label{thm4:eqn2}
\Vert B_{k,\chi}\Vert_{L^{\infty}}=\sup_{\tau\in\mathbb{H}}\Vert B_{k,\chi}(\tau,\tau)\Vert_{\mathrm{Pet}}=O_{\Gamma}\big(k^{\frac{3}{2}}\big).
\end{align}
The implied constants in the $L^{\infty}$-norm bounds~\eqref{thm4:eqn1} and~\eqref{thm4:eqn2} depend only on the Fuchsian subgroup $\Gamma$. 
\begin{proof}
When $\Gamma$ is cocompact without elliptic elements, the claimed $L^{\infty}$-norm bound~\eqref{thm4:eqn1} follows directly from Proposition
\ref{prop1}. \\
Let next $\Gamma$ be a cofinite Fuchsian subgroup. From the proof of Theorem~6.1 in~\cite{fjk2}, it follows that 
\begin{align}
\label{thm4-eqn1}
\sup_{\tau\in\mathbb{H}}\Vert B_{k,\chi}(\tau,\tau)\Vert_{\mathrm{Pet}}=\sup_{\substack{\tau\in\partial\mathcal{F}_{Y}\\Y=k/(4\pi)}}\Vert B_{k,\chi}(\tau,
\tau)\Vert_{\mathrm{Pet}},
\end{align}
where $\partial\mathcal{F}_{Y}$ denotes the boundary of the truncated fundamental domain $\mathcal{F}_{Y}$ defined in \eqref{calfy}. Combining
Proposition~\ref{prop3} with~\eqref{thm4-eqn1} and employing the fact that
\begin{align*}
\frac{\Gamma\big((k-1)/2\big)}{\Gamma(k/2)}=O\bigg(\frac{1}{\sqrt{k}}\bigg),
\end{align*}
we derive 
\begin{align*}
&\sup_{\tau\in\mathbb{H}}\Vert B_{k,\chi}(\tau,\tau)\Vert_{\mathrm{Pet}} \\
&\qquad\le\frac{k-1}{2\pi}+\frac{3(k-1)}{\pi\cosh^{k-4}(r_{\Gamma}/4)}\bigg(1+\frac{1}{\sinh^{2}(r_{\Gamma}/4)}\bigg)+(k-1)\,C_{\Gamma,\mathrm
{ell}}+ k^{\frac{3}{2}}\,C_{\Gamma,\mathrm{par}}
\end{align*}
for some positive constants $C_{\Gamma,\mathrm{ell}},C_{\Gamma,\mathrm{par}}$, which depend only on the Fuchsian subgroup $\Gamma$. 
This completes the proof of the theorem.
\end{proof}
\end{thm}

\begin{cor}
\label{cor5}
With notations as above, let $\Gamma\subset\mathrm{SL_{2}(\mathbb{R})}$ be a cofinite Fuchsian subgroup. For $k\in\mathbb{R}_{\ge 5}$, let 
$f\in S_{k,\chi}(\Gamma)$ be an $L^{2}$-normalized cusp form. If $\Gamma$ is cocompact without elliptic elements, we have the $L^{\infty}$-norm
bound
\begin{align}
\label{cor5:eqn1}
\Vert f\Vert_{L^{\infty}}=\sup_{\tau\in\mathbb{H}}\Vert f(\tau)\Vert_{\mathrm{Pet}}^{2}=O_{\Gamma}(k).
\end{align}
Moreover, if $\Gamma$ is cofinite, we have the $L^{\infty}$-norm bound
\begin{align}
\label{cor5:eqn2}
\Vert f\Vert_{L^{\infty}}=\sup_{\tau\in\mathbb{H}}\Vert f(\tau)\Vert_{\mathrm{Pet}}^{2}=O_{\Gamma}\big(k^{\frac{3}{2}}\big).
\end{align}
The implied constants in the $L^{\infty}$-norm bounds~\eqref{cor5:eqn1} and~\eqref{cor5:eqn2} depend only on the Fuchsian subgroup~$\Gamma$.  
\begin{proof}
Choose an orthonormal basis $\{f_{1}=f,\ldots,f_{d_{k}}\}$ of $S_{k,\chi}(\Gamma)$ with respect to the Petersson inner product~\eqref{pet-ip-gamma}. 
For $\tau\in\mathbb{H}$, we then have the bound
\begin{align*}
\Vert f(\tau)\Vert_{\mathrm{Pet}}^{2}\le\Vert  B_{k,\chi}(\tau,\tau)\Vert_{\mathrm{Pet}}.
\end{align*}
The proof of the corollary now immediately follows from Theorem~\ref{thm4}.
\end{proof}
\end{cor}

\begin{thm}
\label{thm6}
With notations as above, let $\Gamma_{0}\subset\mathrm{SL_{2}}(\mathbb{R})$ be a fixed cofinite Fuchsian subgroup and let $\Gamma\subseteq
\Gamma_{0}$ be a finite index subgroup of $\Gamma_{0}$. For $k\in\mathbb{R}_{\ge 5}$, let $f\in S_{k,\chi}(\Gamma)$ be an $L^{2}$-normalized
cusp form. If $\Gamma_{0}$ is cocompact without elliptic elements, we have the $L^{\infty}$-norm bound
\begin{align}
\label{thm6:eqn1}
\Vert f\Vert_{L^{\infty}}=\sup_{\tau\in\mathbb{H}}\Vert f(\tau)\Vert_{\mathrm{Pet}}^{2}=O_{\Gamma_{0}}(k).
\end{align}
Moreover, if $\Gamma_{0}$ is cofinite, we have the $L^{\infty}$-norm bound
\begin{align}
\label{thm6:eqn2}
\Vert f\Vert_{L^{\infty}}=\sup_{\tau\in\mathbb{H}}\Vert f(\tau)\Vert_{\mathrm{Pet}}^{2}=O_{\Gamma_{0}}\big(k^{\frac{3}{2}}\big).
\end{align}
The implied constants in the $L^{\infty}$-norm bounds~\eqref{thm6:eqn1} and~\eqref{thm6:eqn2} depend only on the Fuchsian subgroup~$\Gamma_
{0}$.  
\begin{proof}
Choose an orthonormal basis $\{f_{1}=f,\ldots,f_{d_{k}}\}$ of $S_{k,\chi}(\Gamma)$ with respect to the Petersson inner product~\eqref{pet-ip-gamma}. 

Let now $\Gamma_{0}$ be a cocompact Fuchsian subgroup without elliptic elements. From the proof of Proposition~\ref{prop1},we derive the bound
\begin{align*}
\sup_{\tau\in\mathbb{H}}\Vert f(\tau)\Vert_{\mathrm{Pet}}^{2}&\le\sup_{\tau\in\mathbb{H}}\Vert B_{k,\chi}(\tau,\tau)\Vert_{\mathrm{Pet}} \\[1mm]
\notag
&\le\frac{k-1}{4\pi}\sup_{\tau\in\mathbb{H}}\sum_{\gamma\in\Gamma}\frac{1}{\cosh^{k}(\mathrm{dist}_{\mathrm{hyp}}(\tau,\gamma\tau)/2)} \\
&\le\frac{k-1}{4\pi}\sup_{\tau\in\mathbb{H}}\sum_{\gamma\in\Gamma_{0}}\frac{1}{\cosh^{k}(\mathrm{dist}_{\mathrm{hyp}}(\tau,\gamma\tau)/2)}=
O_{\Gamma_{0}}(k),
\end{align*}
which completes the proof of the theorem in the case that $\Gamma_{0}$ is cocompact without elliptic elements. 

Let next $\Gamma_{0}$ be a cofinite Fuchsian subgroup. Given $Y>0$, we recall from~\eqref{calfy} the fundamental domain decomposition
\begin{align*}
\mathcal{F}_{\Gamma}=\mathcal{F}_{Y}\cup\big(\mathcal{F}^{Y}_{1}\cup\ldots\cup\mathcal{F}^{Y}_{s}\big),
\end{align*}
where $\mathcal{F}_{Y}$ is a compact subset of $\mathcal{F}_{\Gamma}$ and the $\mathcal{F}^{Y}_{j}$'s are neighborhoods of the cusps $p_{j}
\in\mathcal{P}_{\Gamma}$ ($j=1,\ldots,s$). Choosing $Y$ large enough, we can assume without loss of generality in the sequel that the neighborhoods 
$\mathcal{F}^{Y}_{j}$ are pairwise disjoint. We now first provide a bound for the pointwise Petersson norm of $f$, when $\tau$ ranges across the 
compact set $\mathcal{F}_{Y}$, and subsequently we compute bounds for the pointwise Petersson norm of $f$, when $\tau$ ranges across the 
neighborhoods $\mathcal{F}^{Y}_{j}$ of the cusps for fixed, large enough $Y$.

Adapting arguments from the proof of Proposition~\ref{prop3}, we obtain the bound
\begin{align}
\notag
&\sup_{\tau\in\mathcal{F}_{Y}}\Vert B_{k,\chi}(\tau,\tau)\Vert_{\mathrm{Pet}}\le\frac{k-1}{4\pi}\sup_{\tau\in \mathcal{F}_{Y}}\sum_{\gamma\in\Gamma}
\frac{1}{\cosh^{k}(\mathrm{dist}_{\mathrm{hyp}}(\tau,\gamma\tau)/2)} \\[1mm]
\notag
&\qquad\,\le\frac{k-1}{4\pi}\sup_{\tau\in \mathcal{F}_{Y}}\sum_{\gamma\in\Gamma_{0}}\frac{1}{\cosh^{k}(\mathrm{dist}_{\mathrm{hyp}}(\tau,\gamma
\tau)/2)} \\[1mm]
\label{thm6-eqn3}
&\qquad\le\frac{k-1}{2\pi}+\frac{3(k-1)}{\pi\cosh^{k-4}(r_{\Gamma_{0},Y}/4)}\bigg(1+\frac{1}{\sinh^{2}(r_{\Gamma_{0},Y}/4)}\bigg)+\frac{k-1}{4\pi}\sum_
{e_{j}\in\mathcal{E}_{\Gamma_{0}}}(m_{e_{j}}-1),
\end{align}
where
\begin{align*}
r_{\Gamma_{0},Y}=\inf\bigg\{\mathrm{dist}_{\mathrm{hyp}}(\tau,\gamma\tau)\,\bigg\vert\,\tau\in\mathcal{F}_{Y},\,\gamma\in\Gamma_{0}\setminus
\bigcup_{e_{j}\in\mathcal{E}_{\Gamma_{0}}}\Gamma_{0,e_{j}}\bigg\}>0.
\end{align*}
From this, we immediately conclude that
\begin{align}
\label{thm6-eqn1}
\sup_{\tau\in\mathcal{F}_{Y}}\Vert f(\tau)\Vert_{\mathrm{Pet}}^{2}\le\sup_{\tau\in\mathcal{F}_{Y}}\Vert B_{k,\chi}(\tau,\tau)\Vert_{\mathrm{Pet}}=O_
{\Gamma_{0},Y}(k),
\end{align}
where the implied constant depends on $\Gamma_{0}$ and the choice of $Y$.

We are left to provide bounds for the pointwise Petersson norm of $f$, when $\tau$ ranges across the neighborhoods $\mathcal{F}^{Y}_{j}$ of the 
cusps. For this, we will have to distinguish between the two cases $Y>\frac{k}{4\pi}$ and $Y<\frac{k}{4\pi}$. Without loss of generality, we can 
assume that $j=1$, when $p_{1}\in\mathcal{P}_{\Gamma}$ is the cusp at infinity for $\Gamma$ lying above the cusp $p$ at infinity for $\Gamma_
{0}$ with ramification index $[\Gamma_{0,p}:\Gamma_{p_{1}}]$ (note that $\Gamma_{0,p}$ denotes the stabilizer subgroup of $p$ in $\Gamma_
{0}$). With the above notations, we obtain the inclusion $\Gamma\setminus\Gamma_{p_{1}}\subseteq\Gamma_{0}\setminus\Gamma_{0,p}$.

We first treat the case $Y>\frac{k}{4\pi}$, which implies that $\mathcal{F}_{1}^{Y}\subset\mathcal{F}_{1}^{k/(4\pi)}$. Arguing as in the proof of 
Theorem~6.1 in~\cite{fjk2}, we deduce, recalling the inclusion $\Gamma\setminus\Gamma_{p_{1}}\subseteq\Gamma_{0}\setminus\Gamma_{0,p}$, 
that
\begin{align}
\notag
&\quad\sup_{\tau\in\mathcal{F}_{1}^{Y}}\Vert B_{k,\chi}(\tau,\tau)\Vert_{\mathrm{Pet}}\le\sup_{\tau\in\mathcal{F}_{1}^{k/(4\pi)}}\Vert B_{k,\chi}
(\tau,\tau)\Vert_{\mathrm{Pet}}\le\sup_{\tau\in\partial\mathcal{F}_{1}^{k/(4\pi)}}\Vert B_{k,\chi}(\tau,\tau)\Vert_{\mathrm{Pet}} \\
\notag
&\le\frac{k-1}{4\pi}\Bigg(\sup_{\tau\in\partial\mathcal{F}_{1}^{k/(4\pi)}}\sum_{\gamma\in\Gamma\setminus\Gamma_{p_{1}}}\frac{1}{\cosh^{k}
(\mathrm{dist}_{\mathrm{hyp}}(\tau,\gamma\tau)/2)}+\sup_{\tau\in\partial\mathcal{F}_{1}^{k/(4\pi)}}\sum_{\gamma\in\Gamma_{p_{1}}}\frac{1}
{\cosh^{k}(\mathrm{dist}_{\mathrm{hyp}}(\tau,\gamma\tau)/2)}\Bigg) \\[1mm]
\notag
&\le\frac{k-1}{4\pi}\Bigg(\sup_{\tau\in\partial\mathcal{F}_{1}^{k/(4\pi)}}\sum_{\gamma\in\Gamma_{0}\setminus\Gamma_{0,p}}\frac{1}{\cosh^{k}
(\mathrm{dist}_{\mathrm{hyp}}(\tau,\gamma\tau)/2)}+\hspace*{-2mm}\sup_{\tau\in\partial\mathcal{F}_{1}^{k/(4\pi)}}\sum_{\gamma\in\Gamma_
{p_{1}}}\frac{1}{\cosh^{k}(\mathrm{dist}_{\mathrm{hyp}}(\tau,\gamma\tau)/2)}\Bigg). \\
\label{thm6-eqn4} 
\end{align}
Arguments similar to the ones used to derive the bound~\eqref{thm6-eqn3}, lead for the first term of~\eqref{thm6-eqn4} to the bound
\begin{align*}
&\frac{k-1}{4\pi}\sup_{\tau\in\partial\mathcal{F}_{1}^{k/(4\pi)}}\sum_{\gamma\in\Gamma_{0}\setminus\Gamma_{0,p}}\frac{1}{\cosh^{k}(\mathrm
{dist}_{\mathrm{hyp}}(\tau,\gamma\tau)/2)} \\[1mm]
&\quad\le\frac{k-1}{2\pi}+\frac{3(k-1)}{\pi\cosh^{k-4}(r_{\Gamma_{0},k/(4\pi)}/4)}\bigg(1+\frac{1}{\sinh^{2}(r_{\Gamma_{0},k/(4\pi)}/4)}\bigg)+\frac
{k-1}{4\pi}\sum_{e_{j}\in\mathcal{E}_{\Gamma_{0}}}(m_{e_{j}}-1),
\end{align*}
where
\begin{align*}
r_{\Gamma_{0},k/(4\pi)}=\inf\bigg\{\mathrm{dist}_{\mathrm{hyp}}(\tau,\gamma\tau)\,\bigg\vert\,\tau\in\partial\mathcal{F}_{1}^{k/(4\pi)},\,\gamma
\in\Gamma_{0}\setminus\bigg(\Gamma_{0,p}\cup\bigcup_{e_{j}\in\mathcal{E}_{\Gamma_{0}}}\Gamma_{0,e_{j}}\bigg)\bigg\}>0.
\end{align*}
Since it is easy to see that
\begin{align*}
\frac{1}{\sinh^{2}(r_{\Gamma_{0},k/(4\pi)})}=O_{\Gamma_{0}}(1),
\end{align*}
we arrive for the first term of~\eqref{thm6-eqn4} at the bound
\begin{align}
\label{thm6-eqn5}
\frac{k-1}{4\pi}\sup_{\tau\in\partial\mathcal{F}_{1}^{k/(4\pi)}}\sum_{\gamma\in\Gamma_{0}\setminus\Gamma_{0,p}}\frac{1}{\cosh^{k}(\mathrm
{dist}_{\mathrm{hyp}}(\tau,\gamma\tau)/2)}=O_{\Gamma_{0}}(k).
\end{align}
Using the same arguments as in the proof of Proposition~\ref{prop3}, we derive for the second term of~\eqref{thm6-eqn4} the bound
\begin{align}
\label{thm6-eqn6}
\frac{k-1}{4\pi}\sup_{\tau\in\partial\mathcal{F}_{1}^{k/(4\pi)}}\sum_{\gamma\in\Gamma_{p_{1}}}\frac{1}{\cosh^{k}(\mathrm{dist}_{\mathrm{hyp}}
(\tau,\gamma\tau)/2)}=O\big(k^{\frac{3}{2}}\big).
\end{align}
By means of~\eqref{thm6-eqn4}, we thus deduce from~\eqref{thm6-eqn5} and~\eqref{thm6-eqn6} in the case $Y>\frac{k}{4\pi}$ the bound
\begin{align}
\label{thm6-eqn7}
\sup_{\tau\in\mathcal{F}_{1}^{Y}}\Vert f(\tau)\Vert_{\mathrm{Pet}}^{2}\le\sup_{\tau\in\mathcal{F}_{1}^{Y}}\Vert B_{k,\chi}(\tau,\tau)\Vert_{\mathrm
{Pet}}=O_{\Gamma_{0}}\big(k^{\frac{3}{2}}\big).
\end{align}
We finally turn to the case $Y<\frac{k}{4\pi}$, which implies that $\mathcal{F}_{1}^{k/(4\pi)}\subset\mathcal{F}_{1}^{Y}$. Again, arguing as in the 
proof of  Theorem~6.1 in~\cite{fjk2}, we deduce after recalling once again the inclusion $\Gamma\setminus\Gamma_{p_{1}}\subseteq\Gamma_
{0}\setminus\Gamma_{0,p}$ that
\begin{align}
\notag
\sup_{\tau\in\mathcal{F}_{1}^{Y}}\Vert B_{k,\chi}(\tau,\tau)\Vert_{\mathrm{Pet}}&\le\sup_{\tau\in\mathcal{F}_{1}^{Y}\setminus\mathcal{F}_{1}^
{k/(4\pi)}}\Vert B_{k,\chi}(\tau,\tau)\Vert_{\mathrm{Pet}} \\[1mm]
\notag
&\le\frac{k-1}{4\pi}\sup_{\tau\in\mathcal{F}_{1}^{Y}\setminus\mathcal{F}_{1}^{k/(4\pi)}}\sum_{\gamma\in\Gamma_{0}\setminus\Gamma_{0,p}}
\frac{1}{\cosh^{k}(\mathrm{dist}_{\mathrm{hyp}}(\tau,\gamma\tau)/2)} \\[1mm]
\label{thm6-eqn8}
&\quad\,+\frac{k-1}{4\pi}\sup_{\tau\in\mathcal{F}_{1}^{Y}\setminus\mathcal{F}_{1}^{k/(4\pi)}}\sum_{\gamma\in\Gamma_{p_{1}}}\frac{1}{\cosh^
{k}(\mathrm{dist}_{\mathrm{hyp}}(\tau,\gamma\tau)/2)}.
\end{align}
As before, we now obtain the bounds
\begin{align}
\notag
&\frac{k-1}{4\pi}\sup_{\tau\in\mathcal{F}_{1}^{Y}\setminus\mathcal{F}_{1}^{k/(4\pi)}}\sum_{\gamma\in\Gamma_{0}\setminus\Gamma_{0,p}}\frac
{1}{\cosh^{k}(\mathrm{dist}_{\mathrm{hyp}}(\tau,\gamma\tau)/2)} \\[1mm]
\notag
&\quad\le\frac{k-1}{2\pi}+\frac{3(k-1)}{\pi\cosh^{k-4}(r'_{\Gamma_{0},Y}/4)}\bigg(1+\frac{1}{\sinh^{2}(r'_{\Gamma_{0},Y}/4)}\bigg)+\frac{k-1}{4\pi}
\sum_{e_{j}\in\mathcal{E}_{\Gamma_{0}}}(m_{e_{j}}-1) \\[1mm]
\label{thm6-eqn10}
&\quad=O_{\Gamma_{0},Y}(k),
\end{align}
noting that
\begin{align*}
r'_{\Gamma_{0},Y}=\inf\bigg\{\mathrm{dist}_{\mathrm{hyp}}(\tau,\gamma\tau)\,\bigg\vert\,\tau\in\mathcal{F}_{1}^{Y},\,\gamma\in\Gamma_{0}
\setminus\bigg(\Gamma_{0,p}\cup\bigcup_{e_{j}\in\mathcal{E}_{\Gamma_{0}}}\Gamma_{0,e_{j}}\bigg)\bigg\}>0,
\end{align*}
as well as the bound
\begin{align}
\label{thm6-eqn11}
\frac{k-1}{4\pi}\sup_{\tau\in\mathcal{F}_{1}^{Y}\setminus\mathcal{F}_{1}^{k/(4\pi)}}\sum_{\gamma\in\Gamma_{p_{1}}}\frac{1}{\cosh^{k}(\mathrm
{dist}_{\mathrm{hyp}}(\tau,\gamma\tau)/2)}=O\big(k^{\frac{3}{2}}\big).
\end{align}
By means of~\eqref{thm6-eqn8}, we thus deduce from~\eqref{thm6-eqn10} and~\eqref{thm6-eqn11} in the case $Y<\frac{k}{4\pi}$ the bound
\begin{align}
\label{thm6-eqn12}
\sup_{\tau\in\mathcal{F}_{1}^{Y}}\Vert f(\tau)\Vert_{\mathrm{Pet}}^{2}\le\sup_{\tau\in\mathcal{F}_{1}^{Y}}\Vert B_{k,\chi}(\tau,\tau)\Vert_{\mathrm
{Pet}}=O_{\Gamma_{0},Y}\big(k^{\frac{3}{2}}\big).
\end{align}
Since $Y$ has been fixed, the claim of the theorem follows from~\eqref{thm6-eqn1},~\eqref{thm6-eqn7}, and~\eqref{thm6-eqn12}.
\end{proof}
\end{thm}

\begin{rem}
\label{rem7}
If $f\in S_{k,\chi}(\Gamma)$ is not a Hecke eigenform, then there is no evidence from the literature to suggest that the estimates~\eqref{cor5:eqn1} 
and~\eqref{cor5:eqn2} can be improved. Thus, the estimates~\eqref{cor5:eqn1} and~\eqref{cor5:eqn2} are expected to be optimal. 
\end{rem}

\section{Sup-norm bounds for Jacobi cusp  forms}
\label{sec-4}

For $k\in\mathbb{Z}_{\ge 5}$ and $m\in\mathbb{Z}_{\ge 1}$, let $\phi\in J_{k,m}^{\mathrm{cusp}}(\Gamma_{0})$ be an $L^{2}$-normalized Jacobi 
cusp form of weight $k$ and index $m$ for the full modular group $\Gamma_{0}=\mathrm{SL}_{2}(\mathbb{Z})$. We next aim at bounding the 
quantity
\begin{align*}
\Vert\phi\Vert_{L^{\infty}}=\sup_{(\tau,z)\in\mathbb{H}\times\mathbb{C}}\Vert\phi(\tau,z)\Vert_{\mathrm{Pet}}.
\end{align*}
By the $\Gamma_{0}\ltimes\mathbb{Z}^{2}$-invariance of the function $\Vert\phi(\tau,z)\Vert_{\mathrm{Pet}}$ on $\mathbb{H}\times\mathbb{C}$,
it suffices to find a bound for the quantity $\Vert\phi(\tau,z)\Vert_{\mathrm{Pet}}$, when $\tau$ is ranging through a fundamental domain $\mathcal
{F}_{\Gamma_{0}}$ for $\Gamma_{0}$ and $z$ is ranging through a fundamental domain of the elliptic curve $\mathbb{C}/(\mathbb{Z}\oplus\tau
\mathbb{Z})$. Without loss of generality, we can assume in the sequel that
\begin{align*}
\mathcal{F}_{\Gamma_{0}}=\{\tau=\xi+i\eta\in\mathbb{H}\,\vert\,-1/2\le\xi\le 1/2,\,\vert\tau\vert\ge 1\},
\end{align*}
and, given $\tau=\xi+i\eta\in\mathcal{F}_{\Gamma_{0}}$, that
\begin{align*}
E_{\tau}=\{z=x+iy\in\mathbb{C}\,\vert\,0\le x\le 1,\,0\le y\le\eta\}.
\end{align*}
Furthermore, in order to arrive at an $L^{\infty}$-bound of the quantity $\Vert\phi(\tau,z)\Vert_{\mathrm{Pet}}$, we will make use of the decomposition
\eqref{eich-zag1}, namely
\begin{align*}
\phi(\tau,z)=\sum_{\mu=0}^{2m-1}h_{\mu}(\tau)\vartheta_{\mu,m}(\tau,z),
\end{align*}
where the functions $h_{\mu}$ are cusp forms of weight $(k-\frac{1}{2})$ with respect to the finite index subgroup $\Gamma_{1}=\Gamma_{0,1}
(4m)$ of $\Gamma_{0}$ and the theta functions $\vartheta_{\mu,m}$ are defined in~\eqref{defn:theta}. 

\begin{lem}
\label{lem-corr}
With notations as above, we have the bound
\begin{align}
\label{bound_thetas}
\sum_{\mu=0}^{2m-1}\Vert\vartheta_{\mu,m}(\tau,z)\Vert^{2}_{\mathrm{Pet}}\le 2m\,\mathrm{Im}(\tau)^{\frac{1}{2}}+O\big(m^{\frac{1}{2}}\big)
\end{align}
for $(\tau,z)\in\mathbb{H}\times\mathbb{C}$.
\begin{proof}
By definition, we have for $(\tau,z)\in\mathbb{H}\times\mathbb{C}$ the equalities
\begin{align*}
\vartheta_{\mu,m}(\tau,z)&=\sum_{n\in\mathbb{Z}}e^{2\pi im\tau(n-\frac{\mu}{2m})^{2}+2\pi iz(2mn-\mu)} \\
&=\sum_{n\in\mathbb{Z}}e^{2\pi imn^{2}\tau-2\pi i\cdot2mn\cdot\frac{\mu\tau}{2m}+\pi i\frac{\mu^{2}\tau}{2m}+2\pi i\cdot 2mn\cdot z-2\pi i\mu z} \\
&=e^{\pi i\frac{\mu^{2}\tau}{2m}-2\pi i\mu z}\,\sum_{n\in\mathbb{Z}}e^{2\pi imn^{2}\tau+2\pi i(z-\frac{\mu\tau}{2m})2mn} \\
&=\vartheta_{0,m}\Big(\tau,z-\frac{\mu\tau}{2m}\Big)\,e^{\pi i\frac{\mu^{2}\tau}{2m}-2\pi i\mu z}.
\end{align*}
Taking absolute values, we get
\begin{align*}
\vert\vartheta_{\mu,m}(\tau,z)\vert=\Big\vert\vartheta_{0,m}\Big(\tau,z-\frac{\mu\tau}{2m}\Big)\Big\vert\,e^{-\pi\frac{\mu^{2}\eta}{2m}+2\pi\mu y}.
\end{align*}
This implies the relation
\begin{align*}
\Vert\vartheta_{\mu,m}(\tau,z)\Vert_{\mathrm{Pet}}&=\vert\vartheta_{\mu,m}(\tau,z)\vert\,\eta^{\frac{1}{4}}\,e^{-\frac{2\pi my^{2}}{\eta}} \\[1mm]
&=\Big\vert\vartheta_{0,m}\Big(\tau,z-\frac{\mu\tau}{2m}\Big)\Big\vert\,\eta^{\frac{1}{4}}\,e^{-\pi\frac{\mu^{2}\eta}{2m}+2\pi\mu y-2\pi m\frac{y^{2}}
{\eta}} \\[1mm]
&=\Big\vert\vartheta_{0,m}\Big(\tau,z-\frac{\mu\tau}{2m}\Big)\Big\vert\,\eta^{\frac{1}{4}}\,e^{-\frac{2\pi m}{\eta}(y-\frac{\mu\eta}{2m})^{2}}.
\end{align*}
From this, we derive for $(\tau,z)\in\mathbb{H}\times\mathbb{C}$ the bound
\begin{align*}
\Vert\vartheta_{\mu,m}(\tau,z)\Vert_{\mathrm{Pet}}&=\Big\vert\vartheta_{0,m}\Big(\tau,z-\frac{\mu\tau}{2m}\Big)\Big\vert\,\eta^{\frac{1}{4}}\,e^
{-\frac{2\pi m}{\eta}(y-\frac{\mu\eta}{2m})^{2}} \\[1mm]
&\le\sum_{n\in\mathbb{Z}}e^{-2\pi mn^{2}\eta-4\pi mn(y-\frac{\mu\eta}{2m})-\frac{2\pi m}{\eta}(y-\frac{\mu\eta}{2m})^{2}}\,\eta^{\frac{1}{4}} \\
&=\sum_{n\in\mathbb{Z}}e^{-\frac{2\pi m}{\eta}(n\eta+y-\frac{\mu\eta}{2m})^{2}}\,\eta^{\frac{1}{4}}.
\end{align*}
By an integral test, we arrive for $(\tau,z)\in\mathbb{H}\times\mathbb{C}$ at the bound
\begin{align*}
\Vert\vartheta_{\mu,m}(\tau,z)\Vert_{\mathrm{Pet}}&\leq\bigg(1+\int_{-\infty}^{\infty}e^{-\frac{2\pi m}{\eta}(n\eta+y-\frac{\mu\eta}{2m})^{2}}\,
\mathrm{d}n\bigg)\eta^{\frac{1}{4}} \\
&=\bigg(1+\frac{1}{\sqrt{2m\eta}}\bigg)\eta^{\frac{1}{4}},
\end{align*}
which leads to
\begin{displaymath}
\sum_{\mu=0}^{2m-1}\Vert\vartheta_{\mu,m}(\tau,z)\Vert^{2}_{\mathrm{Pet}}\le 2m\eta^{\frac{1}{2}}+O\big(m^{\frac{1}{2}}\big).
\end{displaymath}
This completes the proof of the lemma.
\end{proof}
\end{lem}

\begin{lem}
\label{aux-lem}
With notations as above, let $k\in\mathbb{Z}_{\ge 5}$, let $\Gamma\subseteq\Gamma_{0}$ be a subgroup of finite index, whose stabilizer subgroup 
$\Gamma_{i\infty}$ is generated by $\big(\begin{smallmatrix}1&1\\0&1\end{smallmatrix}\big)$, and let $f\in S_{k,\chi}(\Gamma)$ be an $L^{2}
$-normalized cusp form of weight $k$ and character $\chi$ for $\Gamma$. Then, we have the bound
\begin{align}
\label{bound-aux}
\sup_{\tau\in\mathcal{F}_{\Gamma_{0}}}\big(\Vert f(\tau)\Vert_{\mathrm{Pet}}^{2}\,\mathrm{Im}(\tau)^{\frac{1}{2}}\big)=O_{\Gamma_{0}}(k^{2}),
\end{align}
where the implied constant depends on $\Gamma_{0}$.
\begin{proof}
Since the width of the cusp $p_{1}=i\infty$ of $\Gamma$ equals $1$, the neighborhood $\mathcal{F}_{1}^{Y}$ of $p_{1}$ is contained in the 
fundamental domain $\mathcal{F}_{\Gamma_{0}}$, provided that $Y>1$. For a given $L^{2}$-normalized $f\in S_{k,\chi}(\Gamma)$, we then write
\begin{align*}
\Vert f(\tau)\Vert_{\mathrm{Pet}}^{2}\,\mathrm{Im}(\tau)^{\frac{1}{2}}=\bigg\vert\frac{f(\tau)}{e^{2\pi i\tau}}\bigg\vert^{2}\cdot\frac{\eta^{k+1/2}}{e^{4\pi 
\eta}}.
\end{align*}
Now, the function $\vert f(\tau)/e^{2\pi i\tau}\vert^{2}$ is subharmonic and bounded in the neighborhood $\mathcal{F}_{1}^{Y}$ and thus takes its
maximum on the boundary $\partial\mathcal{F}_{1}^{Y}$ by the maximum principle for subharmonic functions. Furthermore, by elementary calculus, 
we find that the function $\eta^{k+1/2}/e^{4\pi\eta}$ attains its maximum at
\begin{align*}
\eta_{0}=\frac{k+1/2}{4\pi},
\end{align*}
and is strictly monotonically decreasing for $\eta>\eta_{0}$. By assuming for the moment that $k>12$ and by choosing
\begin{align*}
Y\coloneqq\eta_{0}=\frac{k+1/2}{4\pi}>1,
\end{align*}
we derive from the arguments used in the proof of Theorem~\ref{thm6}, in particular from the bound~\eqref{thm6-eqn4}, that
\begin{align}
\notag
&\sup_{\tau\in\mathcal{F}_{\Gamma_{0}}}\big(\Vert f(\tau)\Vert_{\mathrm{Pet}}^{2}\,\mathrm{Im}(\tau)^{\frac{1}{2}}\big)=\sup_{\tau\in\mathcal{F}_{1}^
{\eta_{0}}}\big(\Vert f(\tau)\Vert_{\mathrm{Pet}}^{2}\,\eta^{\frac{1}{2}}\big) \\[2mm]
\notag
&\qquad\quad=\sup_{\tau\in\partial\mathcal{F}_{1}^{\eta_{0}}}\big(\Vert f(\tau)\Vert_{\mathrm{Pet}}^{2}\,\eta^{\frac{1}{2}}\big)\le\sup_{\tau\in\partial
\mathcal{F}_{1}^{\eta_{0}}}\big(\Vert B_{k,\chi}(\tau,\tau)\Vert_{\mathrm{Pet}}\,\eta^{\frac{1}{2}}\big) \\
\label{auxlem-eqn1}
&\qquad\quad\le k^{\frac{3}{2}}\bigg(\sup_{\tau\in\partial\mathcal{F}_{1}^{\eta_{0}}}\sum_{\gamma\in\Gamma\setminus\Gamma_{p_{1}}}\frac{1}{\cosh^
{k}(\mathrm{dist}_{\mathrm{hyp}}(\tau,\gamma\tau)/2)}+\sup_{\tau\in\partial\mathcal{F}_{1}^{\eta_{0}}}\sum_{\gamma\in\Gamma_{p_{1}}}\frac{1}{\cosh^
{k}(\mathrm{dist}_{\mathrm{hyp}}(\tau,\gamma\tau)/2)}\bigg).
\end{align}
Finally, again using the arguments given in the proof of Theorem~\ref{thm6}, shows that
\begin{align}
\label{auxlem-eqn2}
&\sup_{\tau\in\partial\mathcal{F}_{1}^{\eta_{0}}}\sum_{\gamma\in\Gamma\setminus\Gamma_{p_{1}}}\frac{1}{\cosh^{k}(\mathrm{dist}_{\mathrm{hyp}}
(\tau,\gamma\tau)/2)}=O_{\Gamma_{0}}(1), \\
\label{auxlem-eqn3}
&\sup_{\tau\in\partial\mathcal{F}_{1}^{\eta_{0}}}\sum_{\gamma\in\Gamma_{p_{1}}}\frac{1}{\cosh^{k}(\mathrm{dist}_{\mathrm{hyp}}(\tau,\gamma\tau)/2)}
=O_{\Gamma_{0}}\big(k^{\frac{1}{2}}\big).
\end{align}
Combining the estimate~\eqref{auxlem-eqn1} with the bounds~\eqref{auxlem-eqn2} and~\eqref{auxlem-eqn3}, establishes the proof of the lemma for
$k>12$. The finitely many remaining cases $k\in\{5,\ldots,12\}$ can now be obviously added to the statement by possibly enlarging the implied constant.
\end{proof}
\end{lem}

\begin{lem}
\label{lem-index}
For the index of the congruence subgroup
\begin{align*}
\Gamma_{1}=\Gamma_{0,1}(4m)=\bigg\{\begin{pmatrix}a&b\\c&d\end{pmatrix}\in\Gamma_{0}\,\bigg\vert\,c\equiv 0,\,d\equiv 1\mod 4m\bigg\}
\end{align*}
in $\Gamma_{0}=\mathrm{SL}_{2}(\mathbb{Z})$, we have
\begin{align*}
[\Gamma_{0}:\Gamma_{1}]=O_{\epsilon}\big(m^{2+\epsilon}\big),
\end{align*}
where the implied constant depends on the choice of $\epsilon>0$.
\begin{proof}
The group $\Gamma_{0,1}(4m)$ fits into the short exact sequence
\begin{align*}
0\longrightarrow\Gamma_{0,1}(4m)\longrightarrow\Gamma_{0}(4m)\longrightarrow(\mathbb{Z}/4m\mathbb{Z})^{\times}\longrightarrow 0,
\end{align*}
where the homomorphism from $\Gamma_{0}(4m)$ to $(\mathbb{Z}/4m\mathbb{Z})^{\times}$ is given by the assignment $\big(\begin{smallmatrix}
a&b\\c&d\end{smallmatrix}\big)\mapsto d$. From this we conclude that
\begin{align*}
[\Gamma_{0}:\Gamma_{1}]&=[\mathrm{SL}_{2}(\mathbb{Z}):\Gamma_{0}(4m)]\cdot[\Gamma_{0}(4m):\Gamma_{0,1}(4m)] \\
&=[\mathrm{SL}_{2}(\mathbb{Z}):\Gamma_{0}(4m)]\cdot\big\vert(\mathbb{Z}/4m\mathbb{Z})^{\times}\big\vert=[\mathrm{SL}_{2}(\mathbb{Z}):\Gamma_
{0}(4m)]\cdot\varphi(m),
\end{align*}
where $\varphi(\cdot)$ denotes Euler's $\varphi$-function. Since we trivially have $\varphi(m)\le m$ and since it is well known that 
\begin{align*}
[\mathrm{SL}_{2}(\mathbb{Z}):\Gamma_{0}(4m)]=O_{\epsilon}\big(m^{1+\epsilon}\big)
\end{align*}
with an implied constant depending on the choice of $\epsilon>0$, the claim follows.
\end{proof}
\end{lem}

\begin{thm}
\label{thm11}
For $k\in\mathbb{Z}_{\ge 5}$ and $m\in\mathbb{Z}_{\ge 1}$, let $\phi\in J_{k,m}^{\mathrm{cusp}}(\Gamma_{0})$ be an $L^{2}$-normalized Jacobi 
cusp form of weight $k$ and index $m$ for the full modular group $\Gamma_{0}=\mathrm{SL}_{2}(\mathbb{Z})$. Then, we have the $L^{\infty}$-norm 
bound
\begin{align}
\label{thm11:eqn}
\Vert\phi\Vert^{2}_{L^{\infty}}=\sup_{(\tau,z)\in\mathbb{H}\times\mathbb{C}}\Vert\phi(\tau,z)\Vert^{2}_{\mathrm{Pet}}=O_{\Gamma_{0}, \epsilon}\big
(k^{2}\,m^{\frac{7}{2}+\epsilon}\big),
\end{align}
where the implied constant depends on $\Gamma_{0}$ and the choice of $\epsilon>0$.
\begin{proof}
Substituting the decomposition~\eqref{eich-zag1} of the Jacobi form $\phi(\tau,z)$ into its pointwise Petersson norm and applying the Cauchy--Schwartz 
inequality, we find the estimate
\begin{align*}
\notag
\Vert\phi(\tau,z)\Vert^{2}_{\mathrm{Pet}}&=\bigg\vert\sum_{\mu=0}^{2m-1}h_{\mu}(\tau)\vartheta_{\mu,m}(\tau,z)\bigg\vert^{2}\,\eta^{k}\,e^{-\frac
{4\pi my^{2}}{\eta}} \\
\notag
&\leq\bigg(\sum_{\mu=0}^{2m-1}\vert h_{\mu}(\tau)\vert^{2}\,\eta^{k-\frac{1}{2}}\bigg)\bigg(\sum_{\mu=0}^{2m-1}\vert\vartheta_{\mu,m}(\tau,z)
\vert^{2}\,\eta^{\frac{1}{2}}\,e^{-\frac{4\pi my^{2}}{\eta}}\bigg) \\[1mm]
\notag
&=\bigg(\sum_{\mu=0}^{2m-1}\Vert h_{\mu}(\tau)\Vert^{2}_{\mathrm{Pet}}\bigg)\bigg(\sum_{\mu=0}^{2m-1}\Vert\vartheta_{\mu,m}(\tau,z)\Vert_
{\mathrm{Pet}}^{2}\bigg).
\end{align*}
Thus, we arrive by means of Lemma~\ref{lem-corr} at the bound
\begin{align*}
\notag
\Vert\phi(\tau,z)\Vert^{2}_{\mathrm{Pet}}&\le\sum_{\mu=0}^{2m-1}\Vert h_{\mu}(\tau)\Vert^{2}_{\mathrm{Pet}}\cdot\sum_{\mu=0}^{2m-1}\Vert
\vartheta_{\mu,m}(\tau,z)\Vert_{\mathrm{Pet}}^{2} \\
&\le 2m\sum_{\mu=0}^{2m-1}\Vert h_{\mu}(\tau)\Vert^{2}_{\mathrm{Pet}}\,\eta^{\frac{1}{2}}+O\bigg(m^{\frac{1}{2}}\,\sum_{\mu=0}^{2m-1}\Vert h_
{\mu}(\tau)\Vert^{2}_{\mathrm{Pet}}\bigg),
\end{align*}
from which we derive, after recalling the comment at the beginning of this section,
\begin{align}
\notag
\Vert\phi\Vert^{2}_{L^{\infty}}&=\sup_{(\tau,z)\in\mathbb{H}\times\mathbb{C}}\Vert\phi(\tau,z)\Vert^{2}_{\mathrm{Pet}}=\sup_{\substack{\tau\in\mathcal
{F}_{\Gamma_{0}}\\z\in E_{\tau}}}\Vert\phi(\tau,z)\Vert^{2}_{\mathrm{Pet}} \\
\label{bounds_supnorms}
&\le 2m\sum_{\mu=0}^{2m-1}\sup_{\tau\in\mathcal{F}_{\Gamma_{0}}}\big(\Vert h_{\mu}(\tau)\Vert^{2}_{\mathrm{Pet}}\,\eta^{\frac{1}{2}}\big)+O\bigg
(m^{\frac{1}{2}}\sum_{\mu=0}^{2m-1}\sup_{\tau\in\mathcal{F}_{\Gamma_{0}}}\Vert h_{\mu}(\tau)\Vert^{2}_{\mathrm{Pet}}\bigg).
\end{align}
In order to be able to apply Lemma~\ref{aux-lem} and Theorem~\ref{thm6} to the two summands in~\eqref{bounds_supnorms}, respectively, we need 
to $L^{2}$-normalize the modular forms $h_{\mu}$ ($\mu=0,\ldots,2m-1$) under consideration. To do so, we observe that formula~\eqref{pet-ip2} in
conjunction with the $L^{2}$-normalization of  $\phi$ gives
\begin{align}
\label{prop8-eqn3}
\Vert\phi\Vert^{2}_{L^{2}}=\frac{1}{\sqrt{4m}\,[\Gamma_{0}:\Gamma_{1}]}\sum_{\mu=0}^{2m-1}\Vert h_{\mu}\Vert^{2}_{L^{2}}=1.
\end{align}
An application of Lemma~\ref{lem-index} to~\eqref{prop8-eqn3} then gives
\begin{align*}
\sum_{\mu=0}^{2m-1}\Vert h_{\mu}\Vert^{2}_{L^{2}}=\sqrt{4m}\,[\Gamma_{0}:\Gamma_{1}]=O_{\epsilon}\big(m^{\frac{5}{2}+\epsilon}\big)
\end{align*}
with an implied constant depending on the choice of $\epsilon>0$. With regard to the first summand in~\eqref{bounds_supnorms}, we thus obtain by 
Lemma~\ref{aux-lem} the bound
\begin{align}
\notag
2m\sum_{\mu=0}^{2m-1}\sup_{\tau\in\mathcal{F}_{\Gamma_{0}}}\big(\Vert h_{\mu}(\tau)\Vert^{2}_{\mathrm{Pet}}\,\eta^{\frac{1}{2}}\big)&=2m\sum_
{\mu=0}^{2m-1}\sup_{\tau\in\mathcal{F}_{\Gamma_{0}}}\bigg(\frac{\Vert h_{\mu}(\tau)\Vert^{2}_{\mathrm{Pet}}}{\Vert h_{\mu}\Vert^{2}_{L^{2}}}\,\eta^
{\frac{1}{2}}\bigg)\Vert h_{\mu}\Vert^{2}_{L^{2}} \\
\label{52}
&=O_{\Gamma_{0}}\big(2mk^{2}\big)\sum_{\mu=0}^{2m-1}\Vert h_{\mu}\Vert^{2}_{L^{2}}=O_{\Gamma_{0},\epsilon}\big(k^{2}\,m^{\frac{7}{2}+\epsilon}
\big)
\end{align}
with an implied constant depending on $\Gamma_{0}$ and the choice of $\epsilon>0$. With regard to the second summand in~\eqref{bounds_supnorms},
we derive from Theorem~\ref{thm6} the bound 
\begin{align}
\notag
m^{\frac{1}{2}}\sum_{\mu=0}^{2m-1}\sup_{\tau\in\mathcal{F}_{\Gamma_{0}}}\Vert h_{\mu}(\tau)\Vert^{2}_{\mathrm{Pet}}&=m^{\frac{1}{2}}\sum_{\mu=0}^
{2m-1}\sup_{\tau\in\mathcal{F}_{\Gamma_{0}}}\bigg(\frac{\Vert h_{\mu}(\tau)\Vert^{2}_{\mathrm{Pet}}}{\Vert h_{\mu}\Vert^{2}_{L^{2}}}\bigg)\Vert h_{\mu}
\Vert^{2}_{L^{2}} \\
\label{53}
&=O_{\Gamma_{0}}\big(m^{\frac{1}{2}}k^{\frac{3}{2}}\big)\sum_{\mu=0}^{2m-1}\Vert h_{\mu}\Vert^{2}_{L^{2}}=O_{\Gamma_{0},\epsilon}\big(k^{\frac{3}{2}}
\,m^{3+\epsilon}\big)
\end{align}
with an implied constant depending on $\Gamma_{0}$ and the choice of $\epsilon>0$. The bounds~\eqref{52} and~\eqref{53} complete the proof of the
theorem.
\end{proof}
\end{thm}

\begin{rem}
Comparing our bound with the one obtained by P.~Anamby and S.~Das  in~\cite{das}, we obtain the same polynomial growth in $k$, while the
polynomial growth of their bound is better by a factor of $m^{\frac{3}{4}}$.
\end{rem}

\section*{Acknowledgements}
The first-named author expresses his gratitude to the Departments of Mathematics of Technical University of Darmstadt and of Humboldt University 
of Berlin for their hospitality, where this article was realized. The first-named author also acknowledges support from the INSPIRE research grant 
\emph{DST/INSPIRE/04/2015/002263} and the SERB Matrics grant \emph{MTR/2018/000636}. The first-named author also would like to thank 
Prof.~D.~S.~Nagaraj for helpful discussions. The second-named author acknowledges support from the DFG Cluster of Excellence MATH+. The
third-named author acknowledges support from the LOEWE research unit \emph{Uniformized Structures in Arithmetic and Geometry}. Lastly, all 
the three authors would like to acknowledge the hospitality of the International Centre for Theoretical Sciences (ICTS) during a visit for participating 
in the program \emph{Algebraic and Analytic Aspects of Automorphic Forms} (Code: ICTS/aforms2019/02), where this project was initiated.  



\end{document}